\documentclass[11.05pt,a4paper]{amsart}
\usepackage{amssymb,amsmath}
\usepackage{color}
\usepackage{latexsym}
\usepackage{amsthm,amsfonts,amssymb,mathrsfs}
\usepackage{rotating}
\usepackage[leqno]{amsmath}
\usepackage{xspace}
\usepackage[all]{xy}
\usepackage{longtable}

\headsep=1cm \numberwithin{equation}{section}
\newtheorem{theorem}{Theorem}[section]
\newtheorem{lemma}[theorem]{Lemma}
\newtheorem{proposition}[theorem]{Proposition}
\newtheorem{corollary}[theorem]{Corollary}

\theoremstyle{definition}
\newtheorem{definition}[theorem]{Definition}
\theoremstyle{remark}

\newtheorem{example}[theorem]{Example}

\newcommand{\Ass}{\operatorname{Ass}}

\newcommand{\sdepth}{\operatorname{sdepth}}
\newcommand{\Min}{\operatorname{Min}}
\newcommand{\Mon}{\operatorname{Mon}}

\newcommand{\G}{\operatorname{G}}

\newcommand{\Ht}{\operatorname{ht}}

\newcommand{\In}{\operatorname{in}}

\newcommand{\reg}{\operatorname{reg}}

\newcommand{\Supp}{\operatorname{Supp}}

\newcommand{\Ann}{\operatorname{Ann}}

\newcommand{\Sym}{\operatorname{Sym}}

\newcommand{\depth}{\operatorname{depth}}

\newcommand{\Max}{\operatorname{Max}}

\newcommand{\fm}{\frak{m}}
\newcommand{\fp}{\frak{p}}

\newenvironment{prf}[1][Proof]{\begin{proof}[\bf #1]}{\end{proof}}
\begin{document}

\author[S. Bandari, K. Divaani-Aazar and A. Soleyman Jahan]{Somayeh Bandari (Tehran),
Kamran Divaani-Aazar (Tehran) and Ali Soleyman Jahan (Sanandaj)}

\title[Pretty cleanness and filter-regular sequences]
{Pretty cleanness and filter-regular sequences}

\address{S. Bandari, Department of Mathematics, Alzahra
University, Vanak, Post Code 19834, Tehran, Iran.}
\email{somayeh.bandari@yahoo.com}

\address{K. Divaani-Aazar, Department of Mathematics, Alzahra
University, Vanak, Post Code 19834, Tehran, Iran-and-School of
Mathematics, Institute for Research in Fundamental Sciences (IPM),
P.O. Box 19395-5746, Tehran, Iran.}
\email{kdivaani@ipm.ir}

\address{A. Soleyman Jahan, Department of Mathematics, University of Kurdistan,
Post Code 66177-15175, Sanandaj, Iran.}
\email{solymanjahan@gmail.com}


\keywords {Almost clean modules;  clean modules; $d$-sequences; filter-regular sequences; pretty
clean modules.\\ The research of the second and third authors are supported by grants from IPM
(no. 90130212 and no. 90130062, respectively).}

\begin{abstract} Let $K$ be a field and $S=K[x_1,\ldots, x_n]$. Let $I$ be a monomial ideal of
$S$ and $u_1,\ldots, u_r$ be monomials in $S$ which form a filter-regular sequence on $S/I$. We
show that $S/I$ is pretty clean if and only if $S/(I,u_1,\ldots, u_r)$ is pretty clean.
\end{abstract}

\maketitle

\section{Introduction}

Let $R$ be a multigraded Noetherian ring and $M$ a finitely generated multigraded $R$-module.
(Here, ``multigraded'' stands for ``$\mathbb{Z}^n$-graded''.) A basic fact in commutative algebra
says that there exists a finite filtration $$\mathcal{F}: 0=M_0\subset M_1\subset \cdots \subset
M_r=M$$ of multigraded submodules of $M$ such that there are multigraded isomorphisms $M_i/M_{i-1}
\cong R/\frak p_i(-a_i)$  for some $a_i\in \mathbb{Z}^n$ and some multigraded prime ideals
$\frak p_i$ of $R$. Such a filtration of $M$ is called a (multigraded) prime filtration. The set
of prime ideals $\frak p_1,\ldots,\frak p_r$ which define the cyclic quotients of $\mathcal{F}$
will be denoted by $\Supp \mathcal{F}$. It is known (and easy to see) that
$\Ass_RM\subseteq \Supp \mathcal{F}\subseteq \Supp_RM.$

Let $\Min M$ denote the set of minimal prime ideals of $\Supp_RM$. Dress \cite{D} called a prime
filtration $\mathcal{F}$ of $M$ {\em clean} if $\Supp \mathcal{F}=\Min M$. Pretty clean filtrations
were defined as a generalization of clean filtrations by Herzog and Popescu \cite{HP}. A prime
filtration $\mathcal{F}$ is called {\em pretty clean} if for all $i<j$ for which
$\frak p_i\subseteq \frak p_j$, it follows that $\frak p_i=\frak p_j$. If $\mathcal{F}$
is a pretty clean filtration of $M$, then $\Supp \mathcal{F}=\Ass_RM$; see
\cite[Corollary 3.4]{HP}. The converse is not true in general as shown by some examples
in \cite{HP} and \cite{S4}.  The prime filtration $\mathcal{F}$ of $M$ is called
{\em almost clean} if $\Supp \mathcal{F}=\Ass_RM$.

The $R$-module $M$ is called {\em clean} (resp. {\em pretty clean} or {\em almost clean})
if it admits a clean (resp. pretty clean or almost clean) filtration. Obviously, cleanness
implies pretty cleanness and pretty cleanness implies almost cleanness.

Throughout, let $K$ be a field and $I$ a monomial ideal of the polynomial ring $S=K[x_1,
\ldots, x_n]$. In this paper, we always consider the ring $S$ with its standard multigrading.
So, an ideal $J$ of $S$ is multigraded if and only if $J$ is a monomial ideal. When $I$ is
square-free, one has $\Ass_SS/I=\Min S/I$, and so the above three concepts coincide for $S/I$.
If $S/I$ is pretty clean, then \cite[Theorem 6.5]{HP} asserts that Stanley's  conjecture holds
for $S/I$; see the paragraph preceding Theorem 3.6 for the statement of this conjecture.

Let $u_1,\ldots, u_r$ be monomials in $S$. If $u_1,\ldots, u_r$ is a regular sequence on $S/I$,
then by \cite[Theorem 2.1]{R} $S/I$ is pretty clean if and only if $S/(I,u_1,\ldots, u_r)$ is
pretty clean. In this paper, we pursuit this line of research not only for regular sequences,
but also for other special types of sequences of monomials.

We show that the assertion of \cite[Theorem 2.1]{R} is also true for cleanness and almost cleanness.
Also, we prove that if $u_1,\ldots, u_r$ is a filter-regular sequence on $S/I$, then $S/I$ is
pretty clean if and only if $S/(I,u_1,\ldots, u_r)$ is pretty clean. Next, we show that if
$u_1,\ldots, u_r$ form a filter-regular sequence on $S/I$, then Stanley's conjecture
is true for $S/I$ if and only if it is true for $S/(I,u_1,\ldots, u_r)$.

Assume that $u_1,\ldots, u_r$ is a minimal set of generators for $I$. We prove that if either
$u_1,\ldots, u_r$ is a $d$-sequence, proper sequence or strong $s$-sequence (with respect to
the reverse lexicographic order), then $S/I$ is pretty clean.

\section{Regular Sequences}

We begin with the following preliminary results.

\begin{lemma}\label{Ass} Let $R$ be a commutative Noetherian ring, $M$ an $R$-module
and $A$ an Artinian submodule of $M$. Then $$\Ass_RM=\Ass_RA\cup \Ass_RM/A.$$
\end{lemma}

\begin{prf} It is well-known that $$\Ass_RA\subseteq \Ass_RM\subseteq \Ass_RA\cup
\Ass_RM/A.$$ On the other hand, \cite[Lemma 2.2]{BSS} yields that
$$\Ass_RM/A\subseteq \Ass_RM\cup \Supp_RA.$$
But $A$ is Artinian, and so $\Supp_RA=\Ass_RA$. This implies our desired equality.
\end{prf}

\begin{lemma}\label{Art} Let $R$ be a multigraded Noetherian ring, $M$ a multigraded
finitely generated $R$-module and $A$ a multigraded Artinian submodule of $M$. If $M/A$
is pretty clean (resp. almost clean), then $M$ is pretty clean (resp. almost clean) too.
\end{lemma}

\begin{prf} Since $A$ is an Artinian $R$-module, one has $$\Min A=\Ass_RA=
\Supp_RA\subseteq \Max R.$$
So obviously, if $M/A$ is pretty clean, then $M$ is pretty clean too. Also, by Lemma
\ref{Ass}, almost cleanness of $M/A$ implies almost cleanness of $M$.
\end{prf}

We denote the maximal monomial ideal $(x_1,\ldots, x_n)$ of the ring $S=K[x_1,\ldots, x_n]$
by $\fm$.  For an $S$-module $M$, $H_{\fm}^i(M)$ denotes $i$th local cohomology module of $M$
with respect to $\fm$. If $M$ is a multigraded finitely generated $S$-module, then $H_{\fm}^i(M)$
is a multigraded Artinian $S$-module for all $i$.

\begin{example} Lemma \ref{Art} is not true for the cleanness. To this end, let $S=K[x,y]$ and
$I=(x^2,xy)$. Set $M:=S/I$ and $A:=H_{\fm}^0(M)$. Clearly $A=(x)/I$, and so $M/A\cong S/(x)$.
It is easy to see that $M/A$ is clean while $M$ is not clean.
\end{example}

\begin{proposition}\label{gamma} Let $M$ be a multigraded finitely generated $S$-module and $A$ a
multigraded Artinian submodule of $M$. Then $M$ is pretty clean if and only if $M/A$ is pretty clean.
\end{proposition}

\begin{prf} In view of Lemma \ref{Art}, it remains to show that if $M$ is pretty clean, then $M/A$
is pretty clean. Let $$\mathcal{F}: 0=M_0\subset M_1\subset\cdots \subset M_r=M$$ be a pretty clean
filtration of $M$. For any $S$-module $N$, let $\ell_S(N)$ denote the length of $N$.  First,
by induction on $t:=\ell_S(H_{\fm}^0(M))$, we show that $M/H_{\fm}^0(M)$ is pretty clean. For $t=0$,
there is nothing to prove. Now, assume that $t>0$ and the claim holds for $t-1$. Then $H_{\fm}^0(M)
\neq 0$, and so ${\fm}\in \Ass_SM=\Supp \mathcal{F}$. Since the filtration $\mathcal{F}$ is pretty
clean and $\Ann_SM_1\subseteq \fm$, it follows that
$M_1\cong S/\fm$, and so $(M_1:_M\fm^{\infty})=H_{\fm}^0(M)$. Then, one
has $$H_{\fm}^0(\frac{M}{M_1})=\frac{M_1:_M\fm^{\infty}}{M_1}=\frac{H_{\fm}^0(M)}{M_1},$$ and so
$$\ell_S(H_{\fm}^0(\frac{M}{M_1}))=\ell_S(H_{\fm}^0(M))-\ell_S(M_1)=t-1.$$
Obviously, $M/M_1$ is pretty clean, and so by the induction hypothesis,
$\frac{M/M_1}{H_{\fm}^0(M/M_1)}$
is pretty clean. But, $$\frac{\frac{M}{M_1}}{H_{\fm}^0(\frac{M}{M_1})}=\frac{\frac{M}{M_1}}
{\frac{H_{\fm}^0(M)}{M_1}}\cong \frac{M}{H_{\fm}^0(M)},$$ and hence $M/H_{\fm}^0(M)$ is pretty clean.

Since $A$ is a multigraded Artinian submodule of $M$, one has $A\subseteq H_{\fm}^0(M)$. From the
first part of the proof, we conclude that $\frac{M/A}{H_{\fm}^0(M)/A}$ is pretty clean. But
$H_{\fm}^0(M)/A$ is a multigraded Artinian submodule of $M/A$, and so Lemma \ref{Art} implies
that $M/A$ is pretty clean.
\end{prf}

In what follows, we recall some needed notation and facts about monomial ideals. For each subset
$H$ of $S$, let $\Mon H$ denote the set of all monomials in $H$. For any monomial ideal $I$ of $S$,
there is a unique minimal generating
set $\G(I)$ of $I$. Clearly, $\G(I)$ is consisting of finitely many monomials and there is no
divisibility among different elements of $\G(I)$. Also for any non-empty
subset $T$ of $\Mon S$, set $\G(T):=\G((T))$. Clearly, $\G((T))$ is a finite subset of $T$.
A monomial ideal of $S$ is irreducible if and only if it is of the form $(x^{a_1}_{i_1},\ldots,
x^{a_t}_{i_t})$, where $a_i\in \mathbb{N}$ for all $i$; see \cite[Corollary 1.3.2]{HH}.
Moreover, $(x^{a_1}_{i_1}, \ldots, x^{a_t}_{i_t})$ is $(x_{i_1},\ldots, x_{i_t})$-primary and
each monomial ideal can be written as a finite intersection of irreducible monomial ideals. Let $I$
be a monomial ideal of $S$ and $\mathcal P: I=\bigcap_{i=1}^rQ_i$ a primary decomposition of
$I$ such that each $Q_i$ is an irreducible monomial ideal of $S$. We use notation
$T_i(\mathcal P)$ for $\G(\Mon(\cap_{j=1}^{i-1}Q_j\setminus Q_i))$.
Notice that $$T_1(\mathcal P)=\G(\Mon(S\setminus Q_1))=\{1\}.$$

For proving our first theorem, we shall need the following lemma.

\begin{lemma}\label{Ali} \cite[Corollary 2.7]{S2} Let $I$ be a monomial ideal of $S$. The
following conditions are equivalent:
\begin{enumerate}
\item[a)]  $S/I$ is  clean (resp. pretty clean or almost clean).
\item[b)]  There exists a  primary decomposition
$\mathcal P: I=\bigcap_{j=1}^rQ_j$ of $I$, where each $Q_j$ is an
irreducible $\frak p_j$-primary monomial ideal, such that
\begin{enumerate}
\item[i)]  $\Ht \frak p_j\leq\Ht \frak p_{j+1}$ for all $j$ and
$\{\frak p_1,\ldots, \frak p_r\}=\Min S/I$,\;
\ \\
(resp. $\Ht \frak p_j\leq \Ht \frak p_{j+1}$ for all $j$ or $\{\frak p_1,\ldots,
\frak p_r\}=\Ass_SS/I$) and
\item[ii)]  $T_j(\mathcal P)$ is a singleton for all $1\leq j\leq r$.
\end{enumerate}
\end{enumerate}
\end{lemma}

Next, we generalize \cite[Theorem 2.1]{R}. It also extends \cite[Corollary 4.10]{STY}.

\begin{theorem}\label{extend} Let $I$ be a monomial ideal of $S$ and $u_1,\ldots,
u_c\in \Mon S$ a regular sequence on $S/I$. Then $S/I$ is clean (resp. pretty clean
or almost clean) if and only if $S/(I,u_1,\ldots, u_c)$ is clean (resp. pretty clean
or almost clean).
\end{theorem}

\begin{prf} By induction on $c$, it is enough to prove the case $c=1$. Let $u\in \Mon S$
be a non zero-divisor on $S/I$. Without loss of generality, we may and do assume that for
some integer $0\leq t<n$, the only variables that divide $u$ are $x_{t+1},\ldots, x_n$. Then $u=\prod_{i=t+1}^nx_i^{a_i}$ for some natural integers $a_{t+1},\ldots, a_n$ and $I=JS$ for some
monomial ideal $J$ of $S':=K[x_1,\ldots, x_t]$.

First, we show that if $S/I$ is  clean (resp. pretty clean or almost clean), then $S/(I,u)$ is
clean (resp. pretty clean or almost clean). Let $\mathcal{P}: I=\cap_{i=1}^rQ_i$ be a primary
decomposition of $I$ which satisfies the condition b) in Lemma \ref{Ali}. Let $1\leq e\leq r$.
Since $$\Ass_SS/I=\{\frak p_1,\ldots, \frak p_r\}$$ and $\Ass_SS/Q_e=\{\frak p_e\}$, it turns out
that $u$ is also a non zero-divisor on $S/Q_e$. Hence $Q_e=q_eS$ for some irreducible monomial
ideal $q_e$ of $S'$. Obviously, $$\mathcal{P'}:(I,u)=\big(\cap_{i=t+1}^n (Q_1,x_i^{a_i})
\big)\cap \big(\cap_{i=t+1}^n (Q_2,x_i^{a_i})\big)\cap \ldots \cap
\big(\cap_{i=t+1}^n (Q_r,x_i^{a_i})\big)$$
is a primary decomposition of $(I,u)$ and each $(Q_i,x_j^{a_j}) $ is an irreducible
$(\frak p_i,x_j)$-primary monomial ideal. We are going to show that the condition b) in
Lemma \ref{Ali} holds for $\mathcal{P'}$. Clearly, $T_1(\mathcal{P'})$ is a singleton.
For each $t+2\leq i\leq n$, we have  $$\G(\Mon(\cap_{j=t+1}^{i-1}(Q_1,x_j^{a_j})
\setminus(Q_1,x_i^{a_i})))=\G(\Mon((Q_1,\prod_{j=t+1}^{i-1}x_j^{a_j})\setminus(Q_1,x_i^{a_i})))
=\{\prod_{j=t+1}^{i-1}x_j^{a_j}\}.$$
Let $2\leq i \leq r$, $t+1\leq h\leq n$ and assume that $T_i(\mathcal P)=\{v\}$. Since
$$((\cap_{j=1}^{i-1}\cap_{k=t+1}^n(Q_j,x_k^{a_k}))\cap(\cap_{l=t+1}^{h-1}(Q_i,x_l^{a_l})))
\setminus(Q_i,x_h^{a_h})=((\cap_{j=1}^{i-1}(Q_j,\prod_{k=t+1}^nx_k^{a_k}))\cap
(Q_i,\prod_{l=t+1}^{h-1} x_l^{a_l}))\setminus (Q_i,x_h^{a_h}),$$  one has
$$\G(\Mon(((\cap_{j=1}^{i-1}\cap_{k=t+1}^n(Q_j,x_k^{a_k}))
\cap(\cap_{l=t+1}^{h-1}(Q_i,x_l^{a_l})))\setminus (Q_i,x_h^{a_h})))
=\{v\prod_{l=t+1}^{h-1}x_l^{a_l}\}.$$
So, $T_i(\mathcal P')$ is a singleton for all $i$. On the other hand, we can easily
deduce that
$$\Ass_S\frac{S}{(I,u)}=\{(\frak p,x_k)|\frak p\in \Ass_S\frac{S}{I}\   \  \text{and} \
 t+1\leq k\leq n\}\  \  (*),$$
$$\Min \frac{S}{(I,u)}=\{(\frak p,x_k)|\frak p\in \Min \frac{S}{I}\   \  \text{and} \
\ t+1\leq k\leq n\} \   \  (\dag)$$
and $\Ht (\fp,x_k)=\Ht \fp+1 \  \  (\ddag)$ for all $\fp\in\Ass_SS/I$ and all $t+1\leq k\leq n$.
Hence $\mathcal P'$ satisfies the condition b) in Lemma \ref{Ali}.

Conversely, let $S/(I,u)$ be clean (resp. pretty clean or almost clean).  So, $(I,u)$ has a primary
decomposition $\mathcal P$ which satisfies the condition b) in Lemma \ref{Ali}. From $(*)$, we
can conclude that $\mathcal{P}$ has the form
$$\mathcal{P}:(I,u)=(Q_1,x_{j_1}^{h_{j_1}})\cap(Q_2,x_{j_2}^{h_{j_2}})\cap\ldots
\cap(Q_s,x_{j_s}^{h_{j_s}}),$$
where for each $1\leq i\leq s$, $Q_i=q_iS$ for some irreducible monomial ideal $q_i$
of $S',$ $\sqrt{Q_i}\in \Ass_SS/I$ and $j_i\in\{t+1,\ldots, n\}$. It follows
that $I=\cap_{i=1}^sQ_i$ is a primary decomposition of $I$. By deleting unneeded components, we
get a primary decomposition
$$\mathcal{P'}:I=Q_{i_1}\cap Q_{i_2}\cap\ldots \cap Q_{i_l}$$ such that $i_1<i_2<\cdots <i_l$
and for each $1\leq j\leq l$, $\cap_{k<j}Q_{i_{k}}\nsubseteq Q_{i_j}$ and $\cap_{k<j}Q_{i_{k}}=
\cap_{m<i_j}Q_m$. We intend to show that $\mathcal{P'}$ satisfies the condition b) in Lemma
\ref{Ali}. Since $$\Ass_SS/I=\{\sqrt{Q_{i_1}}, \sqrt{Q_{i_2}},\ldots, \sqrt{Q_{i_l}}\},$$ in view
of $(*)$, $(\dag)$ and  $(\ddag)$, we only need to indicate
that each $T_i(\mathcal P')$ is a singleton. Let $1\leq j\leq l$. Since
$\cap_{k<j}Q_{i_{k}}\nsubseteq Q_{i_j}$, it follows that there exists at least a monomial $v$ in
$\G(\cap_{k<j}Q_{i_{k}})\setminus Q_{i_j}$. We claim that $v$ is unique. If there exists a monomial
$w\neq v$ in $\G(\cap_{k<j}Q_{i_{k}})\setminus Q_{i_j}$, then since $\cap_{k<j}Q_{i_{k}}=
\cap_{m<i_j}Q_m$, it turns out that $v$ and $w$ are belonging to $\G(\cap_{m<i_j}Q_m)\setminus Q_{i_j}$.
Denote $i_j$ by $d$. Since $v,w\in S'$, we can conclude that
$v$ and $w$ are belonging to $$\G((Q_1,x_{j_1}^{h_{j_1}})\cap (Q_2,x_{j_2}^{h_{j_2}})\cap \ldots
\cap (Q_{d-1},x_{j_{d-1}}^{h_{j_{d-1}}}))\setminus (Q_d,x_{j_d}^{h_{j_d}}).$$ This contradicts
the assumption that $T_d(\mathcal P)$ is a singleton. Therefore, each $T_i(\mathcal P')$ is a
singleton, as desired.
\end{prf}

As an immediate consequence, we obtain the following result; see \cite[Proposition 2.2]{HSY}.

\begin{corollary}\label{regular} Let $u_1,\ldots, u_t\in \Mon S$ be a regular sequence on
$S$. Then $S/(u_1,\ldots, u_t)$ is clean.
\end{corollary}

\section{Filter-regular Sequences}

\begin{definition} Let $M$ be a multigraded finitely generated $S$-module. A non-unit monomial
$u$ in $S$ is called a {\em filter-regular element} on $M$ if $$u\notin \underset{\fp\in
\Ass_SM-\{\fm\}}\bigcup \fp.$$ A sequence $u_1,\ldots, u_r$ of non-unit monomials in $S$ is called
a {\em filter-regular sequence} on $M$ if for each $1\leq i\leq r,$ $u_i$ is a filter-regular
element on $M/(u_1,\ldots, u_{i-1})M.$
\end{definition}

\begin{lemma}\label{easy} Let $M$ be a multigraded finitely generated $S$-module. An element
$1\neq u\in \Mon S$ is a filter-regular element on $M$ if and only if it is a non zero-divisor
on $M/H^0_\fm(M)$.
\end{lemma}

\begin{prf} Since $H^0_\fm(M)$ is Artinian and $H^0_\fm(\frac{M}{H^0_\fm(M)})=0$,
Lemma 2.1 yields that $$\Ass_S(\frac{M}{H^0_\fm(M)})=\Ass_SM-\{\fm\}.$$
Hence, by definition the claim is immediate.
\end{prf}

\begin{theorem} Let $I$ be a monomial ideal of $S$ and $u_1,\ldots,u_r\in \Mon S$ a
filter-regular sequence on $S/I$. Then $S/I$ is pretty clean if and only if
$S/(I,u_1,\ldots, u_r)$ is pretty clean.
\end{theorem}

\begin{prf} By induction on $r$, it is enough to prove that for a monomial
filter-regular element $u$ on $S/I$,  $S/I$ is pretty clean if and only if
$S/(I,u)$ is pretty clean. For convenience, we set $M:=S/I$. By Proposition 2.4,
$M$ is pretty clean if and only if $M/H_{\fm}^0(M)$ is pretty clean. By Lemma
\ref{easy}, $u$ is a non zero-divisor on $M/H_{\fm}^0(M)$. Hence, in view of the
isomorphism $$\frac{\frac{M}{H_{\fm}^0(M)}}{u(\frac{M}{H_{\fm}^0(M)})}\cong
\frac{M}{uM+H_{\fm}^0(M)},$$ Theorem 2.6 yields that $M/H_{\fm}^0(M)$ is pretty
clean if and only if $\frac{M}{uM+H_{\fm}^0(M)}$ is pretty clean. On the other
hand, as $\frac{uM+H_{\fm}^0(M)}{uM}$ is a multigraded Artinian submodule
of $M/uM,$ by Proposition 2.4 and the isomorphism $$\frac{M}{uM+H_{\fm}^0(M)}\cong
\frac{\frac{M}{uM}}{\frac{uM+H_{\fm}^0(M)}{uM}},$$ it turns out that
$\frac{M}{uM+H_{\fm}^0(M)}$ is pretty clean if and only if $M/uM$ is pretty clean.
Therefore, $M$ is pretty clean if and only if $M/uM$ is pretty clean.
\end{prf}

\begin{corollary} Let monomials $u_1,\ldots, u_r$ be a filter-regular sequence on $S$.
Then $S/(u_1,\ldots, u_r)$ is pretty clean.
\end{corollary}

\begin{lemma}\label{filter Ass} Let $M$ be a multigraded finitely generated $S$-module
and $u_1,\ldots, u_r\in \Mon S$ a filter-regular sequence on $M$. If $\fm\in \Ass_SM$, then
$\fm\in \Ass_S(M/(u_1,\ldots,u_r)M)$.
\end{lemma}

\begin{prf} By induction on $r$, it is enough to prove that if $u$ is a monomial
filter-regular element on $M$ and $\fm\in \Ass_SM$, then $\fm\in \Ass_SM/uM$. Since
$\fm\in \Ass_SM$, there exists $0\neq x\in M$ such that $\fm=0:_Sx$. Then, there exists
a non-negative integer $t$ such that $x\in u^tM\setminus u^{t+1}M$. Hence $x=u^ty$ for
some $y\in M\setminus uM$. Clearly, $0:_Sy\subset S$. Let $\frak p\subset\fm$ be
a prime ideal of $S$ containing $0:_Sy$. Since $u$ is a filter-regular element on $M$
and $\fp\neq \fm$, it follows that $\frac{u}{1}\in S_\frak p$ is $M_{\frak p}$-regular.
Hence $$(0:_Sx)_{\frak p}=0:_{S_\frak p} \frac{u^t}{1}\frac{y}{1}=0:_{S_\frak p} \frac{y}{1}
=(0:_Sy)_\frak p\subseteq \frak pS_\frak p,$$ and so $$(0:_Sx)\subseteq (0:_Sx)_\frak p\cap
S\subseteq \frak p S_\frak p\cap S=\frak p.$$ This is a contradiction, and so $\fm$ is the
unique prime ideal of $S$ containing $(0:_Sy)$.  So, $$\fm=\sqrt{(0:_Sy)}\subseteq
\sqrt{(0:_Sy+uM)}\subset S.$$ Therefore, $\sqrt{(0:_Sy+uM)}=\fm$, and so
$\fm\in \Ass_S M/uM$.
\end{prf}

A decomposition of $S/I$ as direct sum of $K$-vector spaces of the form $\mathcal D: S/I=
\bigoplus_{i=1}^ru_iK[Z_i]$, where $u_i$ is a monomial in $S$ and $Z_i\subseteq \{x_1,\ldots,
x_n\}$, is called a {\em Stanley decomposition} of $S/I$. The number $\sdepth \mathcal D:=
\min\{|Z_i|: i=1,\ldots, r\}$ is called the {\em Stanley depth}
of $\mathcal D$. The {\em Stanley depth} of $S/I$ is defined to be $$\sdepth S/I:=\max \{\sdepth
\mathcal D: \mathcal D \ \text{is a Stanley decomposition of} \  S/I\}.$$
Stanley conjectured \cite{St} that $\depth S/I\leq \sdepth S/I$. This conjecture is known
as Stanley's conjecture. Recently, this conjecture was extensively examined by several
authors; see e.g. \cite{A1}, \cite{A2}, \cite{HP}, \cite{HSY}, \cite{P}, \cite{R}, \cite{S3}
and \cite{S4}. On the other hand, the present third author \cite{S3} conjectured that there
always exists a Stanley decomposition $\mathcal D$ of $S/I$ such that the degree of each $u_i$
is at most $\reg S/I$. We refer to this conjecture as $h$-regularity conjecture. It is known
that for square-free monomial ideals, these two conjectures are equivalent.

\begin{theorem} Let $I$ be a monomial ideal of $S$ and $u_1,\ldots, u_r\in \Mon S$ a
filter-regular sequence on $S/I$. Then Stanley's  conjecture holds for $S/I$ if
and only if it holds for $S/(I,u_1,\ldots, u_r)$.
\end{theorem}

\begin{prf} By induction on $r$, it is enough to prove that if $u$ is a monomial
filter-regular element on $S/I$, then  Stanley's  conjecture holds for $S/I$ if
and only if it holds for $S/(I,u)$. First, assume that $\fm\in \Ass_SS/I$. Then
$\depth S/I=0$ and by Lemma \ref{filter Ass}, $\fm\in \Ass_SS/(I,u)$. So,
$\depth S/(I,u)=0$. Hence the claim is immediate in this case. Now, assume that
$\fm\notin \Ass_SS/I$. Then $u$ is a non zero-divisor on $S/I$, and so
by \cite[Theorem 1.1]{R}, Stanley's  conjecture holds for $S/I$ if and only if
it holds for $S/(I,u)$.
\end{prf}

\section{$d$-Sequences}

\begin{definition} Let $R$ be a commutative Noetherian ring, $M$ a finitely generated
$R$-module and $f_1,\ldots, f_t\in R$.
\begin{enumerate}
\item[i)] $f_1,\ldots, f_t$ is called a {\em $d$-sequence} on $M$ if $f_1,\ldots, f_t$
is a minimal generating set of the ideal $(f_1,\ldots, f_t)$ and
$(f_1,\ldots,f_i)M:_Mf_{i+1}f_k=(f_1,\ldots, f_i)M :_M f_k$
for all $0\leq i <t$ and all  $k\geq i + 1$. A $d$-sequence on $R$ is simply called a
$d$-sequence.
\item[ii)] $f_1,\ldots, f_t$ is called a {\em proper sequence} if
$f_{i+1}H_j(f_1,\ldots, f_i;R)=0$ for all $0\leq i<t$ and all $j >0$. Here $H_j(f_1,\ldots, f_i;R)$
denotes the $j$th Koszul homology of $R$ with respect to $f_1,\ldots, f_i$.
\item[iii)] Let $M=(g_1,\ldots, g_t)$ and $(a_{ij})_{s\times t}$ be a relation matrix of $M$.
Then the symmetric algebra of $M$ is defined by
$\Sym M:=R[y_1,\ldots, y_t]/J,$ where $J=(\sum_{j=1}^t a_{1j}y_j,\ldots, \sum_{j=1}^t a_{sj}y_j)$.
Let $<$ be a monomial order on the monomials in $y_1,\ldots, y_n$ with the property
$y_1<\cdots <y_n$. Set $I_i:=(g_1, \ldots, g_{i-1}):_Sg_i$.
Then $(I_1y_1,\ldots, I_ty_t)\subseteq \In_<J$. The sequence $g_1,\ldots, g_t$ is called an
{\em $s$-sequence} (with respect to $<$) if $(I_1y_1,\ldots, I_ty_t)=\In_<J$.
If in addition $I_1\subseteq\cdots \subseteq I_t$, then $g_1,\ldots, g_t$ is called a
{\em strong $s$-sequence}.
\end{enumerate}
\end{definition}

\begin{definition} Let $I$ be a (not necessarily square-free) monomial ideal of $S$
with $\G(I)=\{u_1,...,u_m\}$. A monomial $u_t$ is called a leaf of $\G(I)$ if $u_t$
is the only element in $\G(I)$ or there exists a $j\neq t$ such that
$\gcd(u_t,u_i)|\gcd(u_t,u_j)$ for all $i\neq t$. In this case, $u_j$ is called a
branch of $u_t$ . We say that $I$ is a monomial ideal of forest type if every non-empty
subset of $\G(I)$ has a leaf.
\end{definition}

\cite[Theorem 1.5]{SZ} yields that if $I$ is a monomial ideal of forest type, then $S/I$
is pretty clean.

\begin{lemma} \label{gcd} Let $u_1,\ldots, u_t$ be a sequence of monomials with
the following properties:
\begin{enumerate}
\item[i)] there is no $i\neq j$ such that $u_i|u_j$; and
\item[ii)] $\gcd(u_i,u_j)|u_k \  \  \text{for all} \   \ 1\leq i<j<k\leq t$.
\end{enumerate}
Then $I=(u_1,\ldots, u_t)$ is of forest type, and so $S/I$ is pretty clean.
\end{lemma}

\begin{prf} For every non-empty subset $A=\{u_{n_1},\ldots, u_{n_s}\}$ of $\{u_1,\ldots, u_t\}$,
we may and do assume that $n_1<n_2<\cdots <n_s$. Then obviously the first element of
$A$ is a leaf and the last element of $A$ is a branch for that leaf. So, $I$ is of forest type.
Then \cite[Theorem 1.5]{SZ} implies that $S/I$ is pretty clean.
\end{prf}

\begin{proposition} Let $I$ be a monomial ideal of $S$ with $\G(I)=\{u_1,\ldots, u_t\}$.
If $u_1,\ldots, u_t$ is a $d$-sequence,  proper sequence or strong $s$-sequence (with
respect to the reverse lexicographic order), then $S/I$ is pretty clean.
\end{proposition}

\begin{prf} By \cite[Corollaries 3.3 and 3.4]{HRT} any $d$-sequence is a strong
$s$-sequence with respect to the reverse lexicographic order and $u_1,\ldots, u_t$ is a proper
sequence if and only if it is a strong $s$-sequence with respect to the reverse lexicographic
order. So, by the hypothesis and \cite[Theorem 3.1]{T},  there is no $i\neq j$ such that
$u_i|u_j$ and $\gcd(u_i,u_j)|u_k \  \  \text{for all}  \  \ 1\leq i<j<k \leq t$. Hence,
by Lemma \ref{gcd}, $S/I$ is pretty clean.
\end{prf}

Let $I$ be a monomial ideal of $S$ and $u$ a monomial which is a $d$-sequence on $S/I$. The
following example shows that it may happen that $S/I$ is pretty clean, but $S/(I,u)$ is not.

\begin{example} Let $I=(x_1x_2,x_2x_3,x_3x_4)$ be a monomial ideal of $S=K[x_1,x_2,x_3,x_4]$.
It is easy to see that $S/I$ is pretty clean and $x_4x_1$ is a $d$-sequence on $S/I$. But, by
\cite[Example 1.11]{S4}, we know that $S/(I,x_4x_1)=S/(x_1x_2,x_2x_3,x_3x_4,x_4x_1)$ is not
pretty clean.
\end{example}

We conclude the paper with the following result.

\begin{corollary} Let $I$ be a monomial ideal of $S$. Assume that either:
\begin{enumerate}
\item[i)] $I$ is generated by a filter-regular sequence; or
\item[ii)] $I$ is  generated by a $d$-sequence.
\end{enumerate}
Then both Stanley's and the $h$-regularity conjectures hold for $S/I$. Also, in each of these cases
$S/I$ is sequentially Cohen-Macaulay and $\depth S/I=\min\{\dim S/\frak p|\frak p\in \Ass_SS/I\}$.
\end{corollary}

\begin{prf} In both cases i) and ii), it follows that $S/I$ is pretty clean; see Corollary 3.4 and
Proposition 4.4.

As $S/I$ is pretty clean, \cite[Theorem 6.5]{HP} asserts that Stanley's  conjecture holds for $S/I$.
In fact, by \cite[Proposition 1.3]{HVZ}, we have $\depth S/I=\sdepth S/I.$  On the other hand, by
\cite[Theorem 4.7]{S3} the $h$-regularity conjecture holds for $S/I$.

Next, as $S/I$ is pretty clean,  \cite[Corollary 4.3]{HP} implies that $S/I$ is sequentially Cohen-Macaulay.
In \cite{S1} this fact is reproved by a different argument and, in addition, it is shown that depth
of $S/I$ is equal to the minimum of the dimension of $S/\frak p$, where $\frak p\in\Ass_SS/I$.
\end{prf}


\end{document}